
\input amstex

\input xy  \xyoption{all}

\documentstyle{amsppt}
   \magnification=1100
   \hsize=6.2truein
   \vsize=9.0truein
   \hoffset 0.1truein
   \parindent=2em


\font\eusm=eusm10                   


\font\eusms=eusm7                       

\font\eusmss=eusm5                      


\newcount\theTime
\newcount\theHour
\newcount\theMinute
\newcount\theMinuteTens
\newcount\theScratch
\theTime=\number\time
\theHour=\theTime
\divide\theHour by 60
\theScratch=\theHour
\multiply\theScratch by 60
\theMinute=\theTime
\advance\theMinute by -\theScratch
\theMinuteTens=\theMinute
\divide\theMinuteTens by 10
\theScratch=\theMinuteTens
\multiply\theScratch by 10
\advance\theMinute by -\theScratch

\def\today{{\number\day\space
 \ifcase\month\or
  January\or February\or March\or April\or May\or June\or
  July\or August\or September\or October\or November\or December\fi
 \space\number\year}}

\define\Aeu{{\mathchoice
     {\text{\eusm A}}
     {\text{\eusm A}}
     {\text{\eusms A}}
     {\text{\eusmss A}}}}

\define\Afr{{\frak A}}

\define\At{{\widetilde A}}

\define\biggnm#1{
  \bigg|\bigg|#1\bigg|\bigg|}

\define\bignm#1{
  \big|\big|#1\big|\big|}

\define\Bt{{\widetilde B}}

\define\Cpx{\bold C}

\define\Deu{{\mathchoice
     {\text{\eusm D}}
     {\text{\eusm D}}
     {\text{\eusms D}}
     {\text{\eusmss D}}}}

\define\eqdef{{\;\overset\text{def}\to=\;}}

\define\Et{{\widetilde E}}

\define\Eto#1{E_{(\to{#1})}}

\define\fpamalg#1{{\dsize\;\operatornamewithlimits*_{#1}\;}}

\define\fpiamalg#1{{\tsize\;({*_{#1}})_{\raise-.5ex\hbox{$\ssize\iota\in I$}}}}

\define\freeprod#1#2{\mathchoice
     {\operatornamewithlimits{\ast}_{#1}^{#2}}
     {\raise.5ex\hbox{$\dsize\operatornamewithlimits{\ast}
      _{#1}^{#2}$}\,}
     {\text{oops!}}{\text{oops!}}}

\define\freeprodi{\mathchoice
     {\operatornamewithlimits{\ast}
      _{\iota\in I}}
     {\raise.5ex\hbox{$\dsize\operatornamewithlimits{\ast}
      _{\sssize\iota\in I}$}\,}
     {\text{oops!}}{\text{oops!}}}

\define\freeprodvni{\mathchoice
      {\operatornamewithlimits{\overline{\ast}}
       _{\iota\in I}}
      {\raise.5ex\hbox{$\dsize\operatornamewithlimits{\overline{\ast}}
       _{\sssize\iota\in I}$}\,}
      {\text{oops!}}{\text{oops!}}}

\define\Ft{{\tilde F}}

\define\GNS{{\text{\rm GNS}}}

\define\Hil{{\mathchoice
     {\text{\eusm H}}
     {\text{\eusm H}}
     {\text{\eusms H}}
     {\text{\eusmss H}}}}


\define\Hilto#1{\Hil_{(\to{#1})}}

\define\induce{\,\downharpoonright}

\define\lambdat{{\widetilde{\lambda}}}

\define\ld#1{{\hbox{..}(#1)\hbox{..}}}

\define\Leu{{\mathchoice
     {\text{\eusm L}}
     {\text{\eusm L}}
     {\text{\eusms L}}
     {\text{\eusmss L}}}}

\define\lrnm#1{\left\|#1\right\|}

\define\Meu{{\mathchoice
     {\text{\eusm M}}
     {\text{\eusm M}}
     {\text{\eusms M}}
     {\text{\eusmss M}}}}

\define\nm#1{\|#1\|}

\define\Nats{\Naturals}

\define\Naturals{{\bold N}}

\define\otdt{\otimes\cdots\otimes}

\define\otdts#1{\otimes_{#1}\cdots\otimes_{#1}}

\define\otimesdd{{\,{\overset{..}\to\otimes}}\,}

\define\oup{^{\text{\rm o}}}

\define\owedge{{
     \operatorname{\raise.5ex\hbox{\text{$
     \ssize{\,\bigcirc\llap{$\ssize\wedge\,$}\,}$}}}}}

\define\owedgeo#1{{
     \underset{\raise.5ex\hbox
     {\text{$\ssize#1$}}}\to\owedge}}

\define\phit{{\tilde\phi}}

\define\pit{{\widetilde\pi}}

\define\Pto#1{{P_{(\to{#1})}}}


\define\pup#1#2{{{\vphantom{#2}}^{#1}\!{#2}}\vphantom{#2}}

\define\QED{\newline
            \line{$\hfill$\qed}\enddemo}

\define\restrict{\lower .3ex
     \hbox{\text{$|$}}}

\define\sigmabar{{\overline\sigma}}

\define\sigmat{{\widetilde\sigma}}

\define\smd#1#2{\underset{#2}\to{#1}}

\define\smdb#1#2{\undersetbrace{#2}\to{#1}}

\define\smdbp#1#2#3{\overset{#3}\to
     {\smd{#1}{#2}}}

\define\smdbpb#1#2#3{\oversetbrace{#3}\to
     {\smdb{#1}{#2}}}

\define\smdp#1#2#3{\overset{#3}\to
     {\smd{#1}{#2}}}

\define\smdpb#1#2#3{\oversetbrace{#3}\to
     {\smd{#1}{#2}}}

\define\smp#1#2{\overset{#2}\to
     {#1}}

\define\thetabar{{\overline\theta}}

\define\tocdots
  {\leaders\hbox to 1em{\hss.\hss}\hfill}

\define\Veu{{\mathchoice
     {\text{\eusm V}}
     {\text{\eusm V}}
     {\text{\eusms V}}
     {\text{\eusmss V}}}}

\define\Weu{{\mathchoice
     {\text{\eusm W}}
     {\text{\eusm W}}
     {\text{\eusms W}}
     {\text{\eusmss W}}}}

\define\Weut{{\widetilde\Weu}}

\define\xit{{\tilde\xi}}


  \newcount\mycitestyle \mycitestyle=1 

  \newcount\bibno \bibno=0
  \def\newbib#1{\advance\bibno by 1 \edef#1{\number\bibno}}
  \ifnum\mycitestyle=1 \def\cite#1{{\rm[\bf #1\rm]}} \fi
  \def\scite#1#2{{\rm[\bf #1\rm, #2]}}


  \newcount\ignorsec \ignorsec=0
  \def\notasec{\ignorsec=1}

  \newcount\secno \secno=0
  \def\newsec#1{\procno=0 \subsecno=0 \ignorsec=0
    \advance\secno by 1 \edef#1{\number\secno}
    \edef\currentsec{\number\secno}}

  \newcount\subsecno
  \def\newsubsec#1{\procno=0 \advance\subsecno by 1
    \edef\currentsec{\number\secno.\number\subsecno}
     \edef#1{\currentsec}}

  \newcount\appendixno \appendixno=0
  \def\newappendix#1{\procno=0 \ignorsec=0 \advance\appendixno by 1
    \ifnum\appendixno=1 \edef\appendixalpha{\hbox{A}}
      \else \ifnum\appendixno=2 \edef\appendixalpha{\hbox{B}} \fi
      \else \ifnum\appendixno=3 \edef\appendixalpha{\hbox{C}} \fi
      \else \ifnum\appendixno=4 \edef\appendixalpha{\hbox{D}} \fi
      \else \ifnum\appendixno=5 \edef\appendixalpha{\hbox{E}} \fi
      \else \ifnum\appendixno=6 \edef\appendixalpha{\hbox{F}} \fi
    \fi
    \edef#1{\appendixalpha}
    \edef\currentsec{\appendixalpha}}

  \newcount\procno \procno=0
  \def\newproc#1{\advance\procno by 1
   \ifnum\ignorsec=0 \edef#1{\currentsec.\number\procno}
                     \edef\currentproc{\currentsec.\number\procno}
   \else \edef#1{\number\procno}
         \edef\currentproc{\number\procno}
   \fi}

  \newcount\subprocno \subprocno=0
  \def\newsubproc#1{\advance\subprocno by 1
   \ifnum\subprocno=1 \edef#1{\currentproc a} \fi
   \ifnum\subprocno=2 \edef#1{\currentproc b} \fi
   \ifnum\subprocno=3 \edef#1{\currentproc c} \fi
   \ifnum\subprocno=4 \edef#1{\currentproc d} \fi
   \ifnum\subprocno=5 \edef#1{\currentproc e} \fi
   \ifnum\subprocno=6 \edef#1{\currentproc f} \fi
   \ifnum\subprocno=7 \edef#1{\currentproc g} \fi
   \ifnum\subprocno=8 \edef#1{\currentproc h} \fi
   \ifnum\subprocno=9 \edef#1{\currentproc i} \fi
   \ifnum\subprocno>9 \edef#1{TOO MANY SUBPROCS} \fi
  }

  \newcount\tagno \tagno=0
  \def\newtag#1{\advance\tagno by 1 \edef#1{\number\tagno}}



\notasec
  \newproc{\QnFalseProp}
  \newproc{\UnivProp}
\newsec{\Embeddings}
  \newproc{\CondExpA}
  \newproc{\GNStensone}
   \newtag{\LtwoPhi}
   \newtag{\LtwoVsum}
   \newtag{\Wigoes}
  \newproc{\EmbThm}
   \newtag{\kappas}
   \newtag{\lamsigmanms}
   \newtag{\restrictednms}
  \newproc{\nogo}
\newsec{\cpmaps}
  \newproc{\tensops}
   \newtag{\pipit}
  \newproc{\fpcpuBB}
   \newtag{\thetathetai}
   \newtag{\thetaan}
   \newtag{\afrat}
   \newtag{\thetaanxi}
   \newtag{\Plthetaan}
\newsec{\vNalgs}
  \newproc{\vNamFP}
  \newproc{\vNamFPrep}
  \newproc{\normality}
   \newtag{\evxev}
  \newproc{\NormalInduced}
  \newproc{\NormalTensor}
  \newproc{\vNCondExpA}
  \newproc{\vNEmbThm}
   \newtag{\vNkappas}
  \newproc{\vNfpcpuBB}
   \newtag{\Advnfp}
   \newtag{\vNthetathetai}
   \newtag{\vNthetaan}

\newbib{\Avitzour}
\newbib{\BocaZZfpcp}
\newbib{\BocaZZfpcpAmalg}
\newbib{\BozejkoLeinertSpeicherZZCond}
\newbib{\BozejkoSpeicherZZPsiIndep}
\newbib{\Ching}
\newbib{\ChodaZZfpcp}
\newbib{\ChodaDykemaZZPIII}
\newbib{\DykemaZZFaithful}
\newbib{\DykemaZZSimplicity}
\newbib{\DykemaZZPII}
\newbib{\DykemaZZExact}
\newbib{\DykemaZZTopEnt}
\newbib{\DykemaHaagerupRordam}
\newbib{\DykemaRordamZZPI}
\newbib{\DykemaRordamZZProj}
\newbib{\KasparovZZStV}
\newbib{\LanceZZHilbertCS}
\newbib{\RieffelZZInduced}
\newbib{\RieffelZZMoritaEqCW}
\newbib{\VoiculescuZZSymmetries}

\topmatter
   \title
     Embeddings of reduced free products of operator algebras
   \endtitle

   \author Etienne F\. Blanchard, Kenneth J\. Dykema
     \thanks K.D.\ was partially supported as an invited researcher funded by
             CNRS of France and by NSF Grant No.\ DMS 0070558.
     \endthanks
   \endauthor

  \date September 11, 2000 \enddate

  \rightheadtext{}

   \leftheadtext{}

   \address Institut de Math\'ematiques de Luminy,
            163, avenue de Luminy, case 907,
            F--13288 Marseille, France
   \endaddress

   \email blanch\@iml.univ-mrs.fr \endemail

   \address Dept.~of Mathematics,
            Texas A\&M University
            College Station TX 77843--3368, USA.
   \endaddress

   \email Ken.Dykema\@math.tamu.edu
          {\it URL:} http://www.math.tamu.edu/\~{\hskip0.1em}Ken.Dykema/
   \endemail

   \abstract
          Given reduced amalgamated free products of C$^*$--algebras
          $(A,\phi)=\freeprodi(A_\iota,\phi_\iota)$ and
          $(D,\psi)=\freeprodi(D_\iota,\psi_\iota)$, an embedding
          $A\hookrightarrow D$ is shown to exist assuming there are
          conditional--expectation--preserving  embeddings
          $A_\iota\hookrightarrow D_\iota$.
          This result is extended to show the existence of the reduced
          amalgamated free product of certain classes of unital completely
          positive maps.
          Finally, analogues of the above mentioned
          results are proved for amagamated free products of von 
Neumann algebras.
   \endabstract

\endtopmatter

\document \TagsOnRight \baselineskip=5pt

\heading Introduction. \endheading

The reduced free product construction, and more generally the reduced 
amalgamated free
product construction for C$^*$--algebras, introduced independently
by Voiculescu~\cite{\VoiculescuZZSymmetries} and (somewhat less 
generally) Avitzour~\cite{\Avitzour},
has received much recent attention.

It is natural to ask:
to what extent does the reduced free product construction satisfy
universal property, analogous to those for the free product of groups or the
full free product of C$^*$--algebras?
Since the reduced free product of C$^*$--algebras frequently gives rise to
simple C$^*$--algebras,
(see~\cite{\Avitzour}, \cite{\DykemaRordamZZPI}, \cite{\DykemaZZSimplicity},
\cite{\DykemaZZPII} and~\cite{\ChodaDykemaZZPIII}), it is clear that any
universal property for the reduced free product should be quite a bit more
restrictive in character than for the full free product;
however, at first glance the following question still seems reasonable.
\proclaim{Question \QnFalseProp}
If
$$ (A,\phi)=\freeprodi(A_\iota,\phi_\iota) $$
is a reduced free product of C$^*$--algebras, where the $\phi_\iota$ are states
on the unital C$^*$--algebras $A_\iota$ having faithful GNS representations,
and if $D$ is a unital C$^*$--algebra with a state $\psi$ and with unital
$*$--homomorphisms $\kappa_\iota:A_\iota\to D$ such that
\roster
\item"(i)" $\psi\circ\kappa_\iota=\phi_\iota$ for every $\iota\in I$,
\item"(ii)" the family $\bigl(\kappa_\iota(A_\iota)\bigr)_{\iota\in I}$ is free
with respect to $\psi$,
\endroster
does it follow that there is a $*$--homomorphism $\kappa:A\to D$ such that,
denoting by $\alpha_\iota:A_\iota\to A$ the injective $*$--homomorphisms
arising from the free product construction,
$\kappa\circ\alpha_\iota=\kappa_\iota$ for every $\iota\in I$?
(Note that $\kappa$ would necessarily be injective.)
\endproclaim

Note that the homomorphism $\kappa$ exists if and only if the GNS
representation $\pi_\psi:D\to\Leu(L^2(D,\psi))$ of $\psi$
is faithful when restricted to the subalgebra of $D$ generated by 
$\bigcup_{\iota\in I}\kappa_\iota(A_\iota)$
As observed in~\scite{\DykemaRordamZZProj}{1.3}, the answer to
Question~\QnFalseProp{} is ``yes'' if the state $\psi$ on $D$ is assumed to be
faithful, (and a similar result holds in the amalgamated case).
However, in general the answer is ``no'', as was shown by the elementary
example~\scite{\DykemaRordamZZProj}{1.4}, (see also the erratum 
to~\cite{\DykemaRordamZZProj}).

The main result of this paper is an embedding result (Theorem~\EmbThm) implying
that the $*$--homomorphism $\pi$ in Question~\QnFalseProp{} does 
exist provided that the free subalgebras
$\pi_\iota(A_\iota)$ lie in free subalgebras of $D$ that taken 
together generate $D$.
Namely, we have the following property.
\proclaim{Property \UnivProp}
Let $I$ be a set and for every $\iota\in I$ let $A_\iota\subseteq D_\iota$ be
a unital inclusion of unital C$^*$--algebras.
Suppose $\phi_\iota$ is a state on $D_\iota$ such that the GNS 
representations of
$\phi_\iota$ and of the restriction 
$\phi_\iota{\restriction}_{A_\iota}$ are faithful.
Consider the reduced free products of C$^*$--algebras,
$$ \aligned
(D,\phi)&=\freeprodi(D_\iota,\phi_\iota) \\
(A,\psi)&=\freeprodi(A_\iota,\phi_\iota{\restriction}_{A_\iota}).
\endaligned $$
Then there is a $*$--homomorphism $\pi:A\to D$ such that for every 
$\iota\in I$ the
diagram
$$ \matrix \format\c&\quad\c\quad&\l\\
D_\iota & \hookrightarrow & D          \\
\cup    &                 & \,\uparrow\pi \\
A_\iota & \hookrightarrow & A
\endmatrix $$
commutes, where the horizontal arrows are the inclusions
arising from the free product construction.
\endproclaim
This property was previously known under the additional assumption that every
$\psi_\iota$ is faithful, which by~\cite{\DykemaZZFaithful} implies that $\psi$
is faithful on $D$;
as noted after Question~\QnFalseProp, this in turn implies the existence of
$\pi$.
Theorem~\EmbThm{} actually proves more generally a version of
Property~\UnivProp{} for reduced amalgamated free products of C$^*$--algebras.
Such an embedding result is frequently useful for understanding reduced free
product C$^*$--algebras;
it has been used in~\cite{\DykemaZZTopEnt} and several times
in~\cite{\DykemaZZExact}.

We should point out that M\.~Choda has in~\cite{\ChodaZZfpcp} stated a theorem
about reduced free products of completely positive maps which is more general
than Property~\UnivProp.
However, her proof is incomplete, as it implicitly uses the full generality of
Property~\UnivProp{} without justifying its validity.

In~\S\Embeddings, the main theorem about embeddings of reduced
amalgamated free products of C$^*$--algebras is proved.
In~\S\cpmaps, Choda's argument proving the existence of reduced free products
of state--preserving completely positive maps is generalized to prove existence
of reduced amalgamated free products of certain sorts of completely positive
maps.
In ~\S\vNalgs, we consider the reduced free product with amalgamation of von
Neumann algebras and prove analogues of the results in~\S\Embeddings{}
and~\S\cpmaps{} for von Neumann algebras.

\proclaim{Acknowledgements.}\rm  Much of this work was done while 
K.D. was visiting the
Institute of Mathematics of Luminy.
He would like to thank the members of the Institute for their
hospitality and for providing stimulating atmosphere during his visit.
The authors would like to thank the referee for helpful comments.
\endproclaim

\heading \S\Embeddings.  Embeddings. \endheading

In this section we prove the main embedding result.
We use the same notation as in~\cite{\DykemaZZExact} for the reduced 
amalgamated free product
construction.

In the following lemma, with the reduced amalgamated free product of
C$^*$--algebras $(A,\phi)=\freeprodi(A_\iota,\phi_\iota)$ we view each
$A_\iota$ as a C$^*$--subalgebra of $A$ via the canonical embedding arising
from the free product construction.

\proclaim{Lemma \CondExpA}
Let $B$ be a unital C$^*$--algebra, let $I$ be a set and for every $\iota\in I$
let $A_\iota$ be a unital C$^*$--algebra containing a copy of $B$ as a unital
C$^*$--subalgebra and having a conditional expectation
$\phi_\iota:A_\iota\to B$ whose GNS representation is faithful.
Let
$$ (A,\phi)=\freeprodi(A_\iota,\phi_\iota) $$
be the reduced amalgamated free product.
Then for every $\iota_0\in I$ there is a conditional expectation
$\Phi_{\iota_0}:A\to A_{\iota_0}$ such that
$\Phi_{\iota_0}{\restriction}_{A_\iota}=\phi_\iota$ for every
$\iota\in I\backslash\{\iota_0\}$ and $\Phi_{\iota_0}(a_1a_2\cdots a_n)=0$
whenever $n\ge2$ and $a_j\in A_{\iota_j}\cap\ker\phi_{\iota_j}$ with
$\iota_1\ne\iota_2,\ldots,\iota_{n-1}\ne\iota_n$.
\endproclaim
\demo{Proof}
We let $E_\iota=L^2(A_\iota,\phi_\iota)$,
$\xi_\iota=\widehat{1_{A_\iota}}\in E_\iota$,
$E_\iota=\xi_\iota B\oplus E_\iota\oup$.
Then $A$ acts (by definition) on the Hilbert
$B$--module
$$ E=\xi B\oplus\bigoplus\Sb n\ge1\\
  \iota_1,\ldots,\iota_n\in I\\
  \iota_1\ne\iota_2,\iota_2\ne\iota_3,\ldots,\iota_{n-1}\ne\iota_n \endSb
  E_{\iota_1}\oup\otimes_B E_{\iota_2}\oup\otdts BE_{\iota_n}\oup. $$
Identify the submodule $\xi B\oplus E_{\iota_0}\oup$ of $E$ with
the Hilbert $B$--module $E_{\iota_0}$ and let $Q_{\iota_0}:E\to E_{\iota_0}$ be
the projection.
Then $\Phi_{\iota_0}(x)=Q_{\iota_0}xQ_{\iota_0}$ has the desired properties.
\QED

\demo{Explication \GNStensone}
Consider the GNS representation
$\bigl(\sigma,L^2(A,\Phi_{\iota_0}),\eta\bigr)=\GNS(A,\Phi_{\iota_0})$
associated with the conditional expectation $\Phi_{\iota_0}:A\to A_{\iota_0}$
found in Lemma~\CondExpA.
Since $A$ is the closed linear span of $B$ and the set of reduced words of the
form $a_1a_2\cdots a_n$ where $a_j\in A_{\iota_j}\cap\ker\phi_{\iota_j}$ and
$\iota_j\ne\iota_{j+1}$, we see that the Hilbert $A_{\iota_0}$--module in the
GNS representation is
$$ L^2(A,\Phi_{\iota_0})=A_{\iota_0}\oplus\bigoplus\Sb n\ge1\\
\iota_1,\ldots,\iota_n\in I\\
\iota_1\ne\iota_2,\ldots,\iota_{n-1}\ne\iota_n \\
\iota_n\ne\iota_0 \endSb
E_{\iota_1}\oup\otdts BE_{\iota_n}\oup\otimes_BA_{\iota_0}. \tag{\LtwoPhi} $$
Moreover, the action $\sigma$ of $A$ on $L^2(A,\Phi_{\iota_0})$ is determined
by its restrictions $\sigma{\restriction}_{A_\iota}$, which are easily
described.

Let $\rho:A_{\iota_0}\to\Leu(\Veu)$ be a unital $*$--homomorphism, for some
Hilbert space $\Veu$.
Then $\sigma\otimes1:A\to\Leu\bigl(L^2(A,\Phi_{\iota_0})\otimes_\rho\Veu\bigr)$
is a $*$--homomorphism;
it is the representation induced, in the sense of
Rieffel~\cite{\RieffelZZInduced}, from $\rho$ up to $A$, with respect to the
conditional expectation $\Phi_{\iota_0}$, and we will denote this induced
representation by $\rho{\induce}^A$.
We have the following explicit description of $\rho{\induce}^A$, obtained by
tensoring~(\LtwoPhi) with $\otimes_\rho\Veu$ on the right.
Writing $\Hil=L^2(A,\Phi_{\iota_0})\otimes_\rho\Veu$ we have
$$ \Hil=\Veu\oplus\bigoplus
\Sb n\ge1 \\
\iota_1,\ldots,\iota_n\in I \\
\iota_1\ne\iota_2,\ldots,\iota_{n-1}\ne\iota_n \\
\iota_n\ne\iota_0 \endSb
E_{\iota_1}\oup\otdts BE_{\iota_n}\oup\otimes_\rho\Veu. \tag{\LtwoVsum} $$
Moreover, the $*$--homomorphism $\sigma\otimes1$ is determined by its
restrictions
$$
\sigma_\iota\eqdef(\sigma\otimes1){\restriction}_{A_\iota}:A_\iota\to\Leu(\Hil),
$$
given as follows.
Consider the Hilbert spaces
$$ \Hil(\iota)=\cases
\dsize(\eta_\iota B\otimes_{\rho{\restriction}_B}\Veu)\oplus\bigoplus
  \Sb n\ge1\\ \iota_1,\ldots,\iota_n\in I\\
  \iota_1\ne\iota_2,\ldots,\iota_{n-1}\ne\iota_n \\
  \iota_n\ne\iota_0,\,\iota_1\ne\iota \endSb
  E_{\iota_1}\oup\otdts BE_{\iota_n}\oup\otimes_\rho\Veu
\quad&\text{if }\iota\ne\iota_0 \\ \vspace{3ex}
\dsize\bigoplus
  \Sb n\ge1\\ \iota_1,\ldots,\iota_n\in I\\
  \iota_1\ne\iota_2,\ldots,\iota_{n-1}\ne\iota_n \\
  \iota_n\ne\iota_0,\,\iota_1\ne\iota_0 \endSb
  E_{\iota_1}\oup\otdts BE_{\iota_n}\oup\otimes_\rho\Veu
&\text{if }\iota=\iota_0,
\endcases $$
where $\eta_\iota B$ is just the Hilbert $B$--module $B$ with identity element
denoted by $\eta_\iota$.
If $\iota\in I\backslash\{\iota_0\}$ let
$$ W_\iota:E_\iota\otimes_B\Hil(\iota)\to\Hil \tag{\Wigoes} $$
be the unitary defined, using the symbol $\otimesdd$ to denote the tensor
product in~(\Wigoes), by
$$ \align
W_\iota:\;&\xi_\iota\otimesdd(\eta_\iota\otimes v)\mapsto v \\
&\zeta\otimesdd(\eta_\iota\otimes v)\mapsto\zeta\otimes v \\
&\xi_\iota\otimesdd(\zeta_1\otdt\zeta_n\otimes v)
  \mapsto\zeta_1\otdt\zeta_n\otimes v \\
&\zeta\otimesdd(\zeta_1\otdt\zeta_n\otimes v)
  \mapsto\zeta\otimes\zeta_1\otdt\zeta_n\otimes v
\endalign $$
whenever $v\in\Veu$, $\zeta\in E_\iota\oup$,
$\zeta_j\in E_{\iota_j}\oup$ and
$\iota\ne\iota_1,\,\iota_1\ne\iota_2,\,\ldots,\,
\iota_{n-1}\ne\iota_n,\,\iota_n\ne\iota_0$.
Then for every $\iota\in I\backslash\{\iota_0\}$ and $a\in A_\iota$, we have
$$ \sigma_\iota(a)=W_\iota(a\otimes1_{\Hil(\iota)})W_\iota^*. $$
Similarly, define the unitary
$$ W_{\iota_0}:\Veu\oplus\bigl(E_\iota\otimes_B\Hil(\iota_0)\bigr)\to\Hil $$
by
$$ \align
W_{\iota_0}:\;&v\oplus0\mapsto v \\
&0\oplus\bigl(\xi_{\iota_0}\otimesdd(\zeta_1\otdt\zeta_n\otimes v)\bigr)
  \mapsto\zeta_1\otdt\zeta_n\otimes v \\
&0\oplus\bigl(\zeta\otimesdd(\zeta_1\otdt\zeta_n\otimes v)\bigr)
  \mapsto\zeta\otimes\zeta_1\otdt\zeta_n\otimes v.
\endalign $$
Then
$$ \sigma_{\iota_0}(a)=
W_{\iota_0}\bigl(\rho(a)\oplus(a\otimes1_{\Hil(\iota_0)})\bigr)W_{\iota_0}^*.
$$
Note that the above description is related to the construction of the
conditionally free product, due to Bo\D zejko and
Speicher~\cite{\BozejkoSpeicherZZPsiIndep}, (see
also~\cite{\BozejkoLeinertSpeicherZZCond}).
\enddemo

\proclaim{Theorem \EmbThm}
Let $B\subseteq\Bt$ be a (not necessarily unital) inclusion of unital 
C$^*$--algebras.
Let $I$ be a set
and for each $\iota\in I$ suppose
$$ \matrix
1_{\At_\iota}\in& \Bt     & \subseteq & \At_\iota \\
                 & \cup &           & \cup  \\
1_{A_\iota}\in  & B       & \subseteq & A_\iota
\endmatrix $$
are inclusions of C$^*$--algebras.
Suppose that $\phit_\iota:\At_\iota\to\Bt$ is a conditional 
expectation such that $\phit_\iota(A_\iota)\subseteq B$
and assume that $\phit_\iota$ and the restriction 
$\phit_\iota{\restriction}_{A_\iota}$ have faithful GNS 
representations,
for all $\iota\in I$.
Let
$$ \aligned
(\At,\phit)&=\freeprodi(\At_\iota,\phit_\iota) \\
(A,\phi)&=\freeprodi(A_\iota,\phit_\iota{\restriction}_{A_\iota})
\endaligned $$
be the reduced amalgamated free products of C$^*$--algebras.
Then there is a unique $*$--homo\-morph\-ism $\kappa:A\to\At$ such that
for every $\iota\in I$ the diagram
$$ \matrix
\At_\iota&\hookrightarrow&\At \\
\cup&&\phantom{\kappa}\uparrow\kappa \\
A_\iota&\hookrightarrow&A
\endmatrix \tag{\kappas} $$
commutes, where the horizontal arrows are the inclusions arising from 
the free product construction.
Moreover, $\kappa$ is necessarily injective.
\endproclaim
\demo{Proof}
Since $A$ is generated by $\bigcup_{\iota\in I}A_\iota$, it is
clear that $\kappa$ will be unique if it exists.
Let $1$ denote the identity element of $\Bt$ and let $p$ be the identity
element of $B$.
If $p\ne1$ then we may replace $B$ by $B+\Cpx(1-p)$ and each $A_\iota$ by
$A_\iota+\Cpx(1-p)$;
hence we may without loss of generality assume that $B$ is a unital
C$^*$--subalgebra of $\Bt$ and thus each $A_\iota$ is a unital
C$^*$--subalgebra of $\At_\iota$.
Let
$$ \align
(\pit_\iota,\Et_\iota,\xit_\iota)&=\GNS(\At_\iota,\phit_\iota), \\
(\pi_\iota,E_\iota,\xi_\iota)&=\GNS(A_\iota,\phi_\iota)
\endalign $$
and
$$ \align
(\Et,\xit)&=\freeprodi(\Et_\iota,\xit_\iota), \\
(E,\xi)&=\freeprodi(E_\iota,\xi_\iota).
\endalign $$
The inclusion $A_\iota\hookrightarrow\At_\iota$ gives an 
inner--product--preserving isometry of
Banach spaces $E_\iota\hookrightarrow\Et_\iota$ sending $\xi_\iota $ to
$\xit_\iota$, and we identify $E_\iota$ with this subspace of $\Et_\iota$
and thereby $E_\iota\oup$ with the subspace of $\Et_\iota\oup$.
This allows canonical identification of the tensor product module
$$ E_{\iota_1}\oup\otdts B E_{\iota_{p-1}}\oup\otimes_B
\Et_{\iota_p}\oup\otdts\Bt\Et_{\iota_n}\oup $$
with a closed subspace of
$\Et_{\iota_1}\oup\otdts\Bt\Et_{\iota_n}\oup$.
Hence, we may and do identify $E$ with the subspace
$$ \xit B\oplus\bigoplus\Sb
n\ge1 \\
\iota_1,\ldots,\iota_n\in I \\
\iota_1\ne\iota_2,\ldots,\iota_{n-1}\ne\iota_n \endSb
E_{\iota_1}\oup\otdts BE_{\iota_n}\oup $$
of $\Et$.
Let $\Afr=\fpiamalg BA_\iota$ be the universal algebraic free product with
amalgamation over $B$.
Let $\sigma:\Afr\to\Leu(E)$, respectively $\sigmat:\Afr\to\Leu(\Et)$, be the
homomorphism extending the homomorphisms $\pi_\iota:A_\iota\to\Leu(E)$,
respectively $\pit_\iota{\restriction}_{A_\iota}:A_\iota\to\Leu(\Et)$,
($\iota\in I$).
In particular, we have $\overline{\sigma(\Afr)}=A$.
In order to show that $\kappa$ exists, it will suffice to show that
$$ \forall x\in\Afr\qquad\nm{\sigmat(x)}\le\nm{\sigma(x)}. $$
Note that the subspace $E$ of $\Et$ is invariant under $\sigmat(\Afr)$ and that
the restriction of $\sigmat(\cdot)$ to $E$ gives $\sigma$.
This implies
$$ \forall x\in\Afr\qquad\nm{\sigmat(x)}\ge\nm{\sigma(x)}, $$
which will in turn imply that $\kappa$ is injective, once it is known to exist.
Let $\tau$ be a faithful representation of $\Bt$ on a Hilbert space $\Weu$.
Consider the Hilbert space $\Et\otimes_\tau\Weu$ and let
$\lambdat:\Leu(\Et)\to\Leu(\Et\otimes_\tau\Weu)$ be the $*$--homomorphism given
by $\lambdat(x)=x\otimes1_\Weu$.
Then $\lambdat$ is faithful, and hence it will suffice to show that
$$ \forall x\in\Afr\qquad\nm{\lambdat\circ\sigmat(x)}\le\nm{\sigma(x)}.
\tag{\lamsigmanms} $$
Our strategy will be to show that $\lambdat\circ\sigmat$ decomposes as a direct
sum of subrepresentations, each of which is of the form
$(\nu{\induce}^A)\circ\sigma$, where $\nu{\induce}^A$ is the
$*$--representation of
$A$ induced from a representation $\nu$ of some $A_\iota$.

Given $n\ge1$ and $\iota_1,\ldots,\iota_n$ with
$\iota_1\ne\iota_2,\ldots,\iota_{n-1}\ne\iota_n$, and given
$p\in\{1,2,\ldots,n\}$, consider the Hilbert space
$$ \align
E_{\iota_1}\oup&\otdts BE_{\iota_{p-1}}\oup\otimes_B K_{\iota_p}
  \otimes_\Bt\Et_{\iota_{p+1}}\oup\otdts\Bt\Et_{\iota_n}\oup\otimes_\tau\Weu=
  \\ \vspace{2ex}
&\eqdef\left(\aligned
&E_{\iota_1}\oup\otdts B
  E_{\iota_{p-1}}\oup\otimes_B\Et_{\iota_p}\oup\otimes_\Bt\Et_{\iota_{p+1}}\oup
  \otdts\Bt\Et_{\iota_n}\oup\otimes_\tau\Weu \\
\ominus\;\;&E_{\iota_1}\oup\otdts B
  E_{\iota_{p-1}}\oup\otimes_BE_{\iota_p}\oup\otimes_B\Et_{\iota_{p+1}}\oup
  \otdts\Bt\Et_{\iota_n}\oup\otimes_\tau\Weu
\endaligned\right)_.
\endalign $$
Heuristically, $K_\iota$ takes the place of $\Et_\iota\ominus E_\iota$, even
when the latter does not make sense.
Then
$$ \Et\otimes_\tau\Weu=(E\otimes_{\tau{\restriction}_B}\Weu)\oplus
\bigoplus\Sb
n\ge1 \\
\iota_1,\ldots,\iota_n\in I \\
\iota_1\ne\iota_2,\ldots,\iota_{n-1}\ne\iota_n \\
p\in\{1,2,\ldots,n\} \endSb
E_{\iota_1}\oup\otdts BE_{\iota_{p-1}}\oup\otimes_BK_{\iota_p}
\otimes_\Bt\Et_{\iota_{p+1}}\oup\otdts\Bt\Et_{\iota_n}\oup\otimes_\tau\Weu. $$
As mentioned earlier, $\sigmat(\Afr)E\subseteq E$ and
$\sigmat(\cdot){\restriction}_E=\sigma(\cdot)$, so
$E\otimes_{\tau{\restriction}_B}\Weu$ is invariant under
$\lambdat\circ\sigmat(\Afr)$, and
$$ \forall x\in\Afr\qquad
\nm{\lambdat\circ\sigmat(x){\restriction}_{E\otimes_\tau\Weu}}=\nm{\sigma(x)}.
$$
Since $\pit_\iota(A_\iota)E_\iota\subseteq E_\iota$, it is not difficult to
check that for every $n\ge1$ and for every $\iota_1,\ldots,\iota_n\in I$ with
$\iota_1\ne\iota_2,\ldots,\iota_{n-1}\ne\iota_n$,
$$ \align
\Weut(\iota_1,\ldots,\iota_n)&\eqdef
\overline{\lambdat\circ\sigmat(\Afr)(K_{\iota_1}\otimes_\Bt
\Et_{\iota_2}\oup\otdts\Bt\Et_{\iota_n}\oup\otimes_\tau\Weu)}= \\ \vspace{2ex}
&\aligned
=&\;(K_{\iota_1}\otimes_\Bt
\Et_{\iota_2}\oup\otdts\Bt\Et_{\iota_n}\oup\otimes_\tau\Weu)\quad\oplus \\
\vspace{1ex}
&\;\oplus\dsize\bigoplus\Sb
q\ge1 \\
\iota_1',\ldots,\iota_q'\in I \\
\iota_1'\ne\iota_2',\ldots,\iota_{q-1}'\ne\iota_q' \\
\iota_q'\ne\iota_1 \endSb
E_{\iota_1'}\oup\otdts BE_{\iota_q'}\oup\otimes_BK_{\iota_1}\otimes_\Bt
\Et_{\iota_2}\oup\otdts\Bt\Et_{\iota_n}\oup\otimes_\tau\Weu.
\endaligned
\endalign $$
Thus
$$ \Et\otimes_\tau\Weu=(E\otimes_{\tau{\restriction}_B}\Weu)\oplus
\bigoplus\Sb
n\ge1 \\
\iota_1,\ldots,\iota_n\in I \\
\iota_1\ne\iota_2,\ldots,\iota_{n-1}\ne\iota_n \endSb
\Weut(\iota_1,\ldots,\iota_n); $$
hence in order to prove the theorem it will suffice to show that for every
choice of $\iota_1,\ldots\iota_n$,
$$ \forall x\in\Afr\qquad
\nm{\lambdat\circ\sigmat(x){\restriction}_{\Weut(\iota_1,\ldots,\iota_n)}}
\le\nm{\sigma(x)}. \tag{\restrictednms} $$
But letting
$\Veu=K_{\iota_1}\otimes_\Bt
\Et_{\iota_2}\oup\otdts\Bt\Et_{\iota_n}\oup\otimes_\tau\Weu$,
letting $\nu:A_{\iota_1}\to\Leu(\Veu)$ be the $*$--homomorphism
$\nu(a)=(\pit_{\iota_1}(a)\otimes
1_{\Et_{\iota_2}\oup\otdts B\Et_{\iota_n}\oup\otimes_\tau\Weu})
{\restriction}_\Veu$,
and appealing to Explication~\GNStensone, it is straightforward to check that
$$ \lambdat\circ\sigmat(\cdot){\restriction}_{\Weut(\iota_1,\ldots,\iota_n)}
=(\nu{\induce}^A)\circ\sigma, $$
where $\nu{\induce}^A$ is the representation of $A$ induced from $\nu$ with
respect to the conditional expectation $\Phi_{\iota_1}:A\to A_{\iota_1}$ found
in Lemma~\CondExpA;
this in turn implies~(\restrictednms).
\QED

\demo{Remark \nogo}
Let us consider for a moment Theorem~\EmbThm{} when the subalgebra $B$ over
which we amalgamate is the scalars, $\Cpx$.
When taking the reduced free product $(A,\phi)=\freeprodi(A_\iota,\phi_\iota)$
of C$^*$--algebras, one usually requires the states $\phi_\iota$ to 
have faithful GNS
representation.
However, one could extend the construction to the case of completely 
general states
$\phi_\iota$;
one then obtains
$$ \freeprodi(A_\iota,\phi_\iota)=
\freeprodi((A_\iota/\ker\pi_\iota),\overset{.}\to{\phi}_\iota), $$
where $\pi_\iota$ is the GNS representation of $\phi_\iota$ and where
$\overset{.}\to{\phi}_\iota$ is the state induced on the quotient
$A_\iota/\ker\pi_\iota$ by $\phi_\iota$.
Thus the canonical $*$--homomorphism $A_\iota\to A$ has the same
kernel as $\pi_\iota$.

As a caveat, we would like to point out that with this relaxed definition of
reduced free product, (allowing $\phi_\iota$ with nonfaithful GNS
representation), the statement of Theorem~\EmbThm{} does not in general hold.
Indeed, if for some $\iota\in I$ $A_\iota=\Cpx\oplus\Cpx$ with $\phi_\iota$
non--faithful, if $\At_\iota=M_2(\Cpx)$ with a unital embedding
$A_\iota\hookrightarrow\At_\iota$ and if
$\phit_\iota$ is a state on $M_2(\Cpx)$ such that
$\phit{\restriction}_{A_\iota}=\phi_\iota$, then $A_\iota\to A$ is
not injective, while $\At_\iota\hookrightarrow\At$ is injective.
This shows that there can be no $*$--homomorphism $\kappa:A\to\At$
making the diagram~(\kappas) commute.
However, there is no problem allowing the $\phit_\iota$ to have nonfaithful GNS
representations, as long as the restrictions $\phi_\iota$ are taken with
faithful GNS representations.
\enddemo

\heading \S\cpmaps.  Completely positive maps. \endheading

M\. Choda~\cite{\ChodaZZfpcp} gave an argument which, when combined with an
embedding result like Property~\UnivProp, proves that if
$\theta_\iota:A_\iota\to D_\iota$ is a unital completely positive map between
unital C$^*$--algebras for every $\iota\in I$, if $\phi_\iota$ and $\psi_\iota$
are states on $A_\iota$ and respectively $D_\iota$, each having faithful GNS
representation, and if $\psi_\iota\circ\theta_\iota=\phi_\iota$ then letting
$$ \align
(A,\phi)&=\freeprodi(A_\iota,\phi_\iota) \\
(D,\psi)&=\freeprodi(D_\iota,\psi_\iota) \endalign $$
be the reduced free products of C$^*$--algebras,
there is a unital completely positive map $\theta:A\to D$ such that
$\theta{\restriction}_{A_\iota}=\theta_\iota$ for every $\iota\in I$,
and such that
$\theta(a_1a_2\cdots a_n)=\theta(a_1)\theta(a_2)\cdots\theta(a_n)$ whenever
$a_j\in A_{\iota_j}\cap\ker\phi_{\iota_j}$ for some $\iota_j\in I$ with
$\iota_1\ne\iota_2,\ldots,\iota_{n-1}\ne\iota_n$.

In this section, we generalize this argument of Choda's to the case of reduced
amalgamated free products of C$^*$--algebras.
The generalization consists of, in essence, replacing Stinespring's dilation
theorem for completely positive maps into bounded operators on a Hilbert space
by Kasparov's generalization~\cite{\KasparovZZStV} to the case of completely
positive maps into the algebra of bounded adjointable operators on a Hilbert
$B$--module (see alternatively the book~\cite{\LanceZZHilbertCS}).
We would like to point out that Theorem~\fpcpuBB{} is quite similar in
appearance to analogous results of F\.~Boca~\cite{\BocaZZfpcp},
\cite{\BocaZZfpcpAmalg} about completely positive maps on universal amalgamated
free products of C$^*$--algebras.
However, the universal and reduced free products of C$^*$--algebras are quite
different in character, and we do not believe that Boca's results can be used
directly to prove Theorem~\fpcpuBB.

\proclaim{Lemma \tensops}
Let $A$ and $B$ be C$^*$--algebras, let $E$ and $\Et$ be Hilbert $A$--modules,
let $F$ and $\Ft$ be Hilbert $B$--modules and let $v\in\Leu(E,\Et)$,
$w\in\Leu(F,\Ft)$.
Suppose $\pi:A\to\Leu(F)$ and $\pit:A\to\Leu(\Ft)$ are $*$--homomorphisms and
suppose that
$$ \forall a\in A\quad\forall\xi\in F\qquad w(\pi(a)\xi)=\pit(a)w(\xi).
\tag{\pipit} $$
Let $E\otimes_\pi F$ and $\Et\otimes_\pit\Ft$ be the interior tensor products.
Then there is an element $v\otimes w\in\Leu(E\otimes_\pi F,\Et\otimes_\pit\Ft)$
such that
$$ \forall\zeta\in E\quad\forall\xi\in F\qquad
(v\otimes w)(\zeta\otimes\xi)=(v\zeta)\otimes(w\xi). $$
If, moreover, $\langle v(\zeta),v(\zeta)\rangle=\langle\zeta,\zeta\rangle$ for
every $\zeta\in E$ and $\langle w(\xi),w(\xi)\rangle=\langle\xi,\xi\rangle$ for
every $\xi\in F$ then
$\langle v\otimes w(\eta),v\otimes w(\eta)\rangle=\langle\eta,\eta\rangle$ for
every $\eta\in E\otimes_\pi F$.
\endproclaim
\demo{Proof}
That $v\otimes w$ is bounded is a standard argument (compare p.~42
of~\cite{\LanceZZHilbertCS}).
Then one sees $(v\otimes w)^*=v^*\otimes w^*$.
The final statement follows using the polarization identity.
\QED

\proclaim{Theorem \fpcpuBB}
Let $B$ be a unital C$^*$--algebra, let $I$ be a set and for every $\iota\in I$
let $A_\iota$ and $D_\iota$ be unital C$^*$--algebras containing copies of $B$
as unital C$^*$--subalgebras and having conditional expectations
$\phi_\iota:A_\iota\to B$, respectively $\psi_\iota:D_\iota\to B$, whose GNS
representations are faithful.
Suppose that for each $\iota\in I$ there is a unital completely positive map
$\theta_\iota:A_\iota\to D_\iota$ that is also a $B$--$B$ bimodule map and
satisfies $\psi_\iota\circ\theta_\iota=\phi_\iota$.
Let
$$ \align
(A,\phi)&=\freeprodi(A_\iota,\phi_\iota) \\
(D,\phi)&=\freeprodi(D_\iota,\psi_\iota)
\endalign $$
be the reduced amalgamated free products of C$^*$--algebras.
Then there is a unital completely positive map $\theta:A\to D$ such 
that for all $\iota\in I$ the diagram
$$\xymatrix{
A_\iota\ar@{^(->}[dd]\ar[rr]|-{\,\Theta_\iota\,}
\ar[dr]|-{\,\phi_\iota\,}&&
D_\iota\ar@{^(->}[dd]\ar[dl]|-{\,\psi_\iota\,}\\
&B&\\
A\ar[ur]|-{\,\phi\,}\ar@{-->}[rr]|-{\,\Theta\,}&&D\ar[ul]|-{\,\psi\,}}$$
commutes, where the vertical inclusions are those arising from the 
free product construction,
and satisfying
$$ \theta(a_1a_2\cdots a_n)=\theta(a_1)\theta(a_2)\cdots\theta(a_n)
\tag{\thetaan} $$
whenever $a_j\in A_{\iota_j}\cap\ker\phi_{\iota_j}$ and
$\iota_1\ne\iota_2$, $\iota_2\ne\iota_3,\,\ldots,\,\iota_{n-1}\ne\iota_n$.
\endproclaim
\demo{Proof}
Note first that the assumptions imply that each $\theta_\iota$ is the identity
map on $B$.
Let
$$ (\pi_\iota,E_\iota,\xi_\iota)=\GNS(D_\iota,\psi_\iota),
\qquad(E,\xi)=\freeprodi(E_\iota,\xi_\iota). $$
(We will usually write simply $d\zeta$ instead of $\pi_\iota(d)\zeta$, when
$d\in D_\iota$ and $\zeta\in E_\iota$.)
Recall that then $E_\iota=\xi_\iota B\oplus E\oup_\iota$, that the action
$\pi_\iota{\restriction}_B$ leaves $E\oup_\iota$ globally invariant, and that
$$ E=\xi B\oplus\bigoplus
  \Sb n\ge1\\\iota_1,\ldots,\iota_n\in I\\
  \iota_1\ne\iota_2,\ldots,\iota_{n-1}\ne\iota_n\endSb
E\oup_{\iota_1}\otdts B E\oup_{\iota_n}. $$
Consider the Hilbert $B$--module
$F_\iota=A_\iota\otimes_{\pi_\iota\circ\theta_\iota}E_\iota$ and the
specified element $\eta_\iota=1\otimes\xi_\iota\in F_\iota$.
Since $\theta_\iota$ restricts to the identity map on $B$, in $F_\iota$ we have
$b\otimes\zeta=1\otimes(b\zeta)$ for every $b\in B$ and $\zeta\in E$.
Consider the unital $*$--homomorphism $\sigma_\iota:A_\iota\to\Leu(F_\iota)$
given by
$$ \forall a',a\in A_\iota\quad\forall\zeta\in E_\iota\qquad
  \sigma_\iota(a')(a\otimes\zeta)=(a'a)\otimes\zeta, $$
(cf.\ page~48 of~\cite{\LanceZZHilbertCS}).
Consider the map $\rho_\iota:\Leu(F_\iota)\to B$ given by
$\rho_\iota(x)=\langle\eta_\iota,x\eta_\iota\rangle$.
If $x\in\Leu(F_\iota)$ and if $b_1,b_2\in B$ then
$\rho_\iota\bigl(\sigma_\iota(b_1)x\sigma_\iota(b_2)\bigr)=b_1\rho_\iota(x)b_2$.
If we use $\sigma_\iota$ to identify $B$ with
$\sigma_\iota(B)\subseteq\Leu(F_\iota)$ then we have that
$\rho_\iota:\Leu(F_\iota)\to B$ is a conditional expectation.
Clearly $L^2(\Leu(F_\iota),\rho_\iota)\cong F_\iota$ and the GNS representation
of $\rho_\iota$ is faithful on $\Leu(F_\iota)$.
We have that $\rho_\iota\circ\sigma_\iota=\phi_\iota$ since
$$ \rho_\iota\circ\sigma_\iota(a)=
  \langle1\otimes\xi_\iota,a\otimes\xi_\iota\rangle=
  \langle\xi_\iota,\theta_\iota(a)\xi_\iota\rangle=
  \psi_\iota\circ\theta_\iota(a)=\phi_\iota(a). \tag{\afrat} $$
Let
$$ (\Meu,\rho)=\freeprodi(\Leu(F_\iota),\rho_\iota) $$
be the reduced amalgamated free product of C$^*$--algebras.
Note that $\Meu\subseteq\Leu(F)$ where
$$ (F,\eta)=\freeprodi(F_\iota,\eta_\iota). $$
By Theorem~\EmbThm{} there is a $*$--homomorphism $\sigma:A\to\Meu$ such that
$\sigma{\restriction}_{A_\iota}=\sigma_\iota$, ($\iota\in I$).

Consider the operator $v_\iota:E_\iota\to F_\iota$ given by
$\zeta\to1\otimes\zeta$, and note that $\langle
v_\iota\zeta,v_\iota\zeta\rangle=\langle\zeta,\zeta\rangle$ for every
$\zeta\in E_\iota$, hence $v_\iota(E\oup_\iota)\subseteq F\oup_\iota$.
A calculation using e.g\. Lemma~5.4 of~\cite{\LanceZZHilbertCS} shows that
there is a bounded operator $F_\iota\to E_\iota$ sending $a\otimes\zeta$ to
$\theta_\iota(a)\zeta$, which is then the adjoint of $v_\iota$.
Hence $v_\iota\in\Leu(E_\iota,F_\iota)$, and clearly $v_\iota^*v_\iota=1$.
Since $\theta_\iota$ is a left $B$--module map, we have for every $b\in B$ and
$\zeta\in E_\iota$ that
$v_\iota(b\zeta)=1\otimes(b\zeta)=b\otimes\zeta=b(v_\iota(\zeta))$.
Therefore, taking direct sums of operators $v_{\iota_1}\otdt v_{\iota_n}$ given
by Lemma~\tensops, we get $v\in\Leu(E,F)$ such that
$\langle v\zeta,v\zeta\rangle=\langle\zeta,\zeta\rangle$ for every
$\zeta\in E$, $v\xi=\eta$ and
$$ v(\zeta_1\otimes\zeta_2\otdt\zeta_n)=
(v_{\iota_1}\zeta_1)\otimes(v_{\iota_2}\zeta_2)\otdt(v_{\iota_n}\zeta_n) $$
whenever $\zeta_j\in E\oup_{\iota_j}$, $\iota_1,\ldots,\iota_n\in I$ and
$\iota_j\ne\iota_{j+1}$.
Let $\theta:A\to\Leu(E)$ be the unital completely positive map
$\theta(x)=v^*\sigma(x)v$.

We will show that~(\thetathetai) commutes and and that~(\thetaan) 
holds, which will furthermore
imply that $\theta(A)\subseteq D$.
In order to show~(\thetathetai), let $w_\iota:E\to E_\iota\otimes_BE(\iota)$
and $y_\iota:F\to F_\iota\otimes_BF(\iota)$ be the unitaries used in the free
product constructions to define the inclusions $A_\iota\hookrightarrow A$ and,
respectively, $\Leu(F_\iota)\hookrightarrow\Meu$.
Note that $v_\iota\bigl(E(\iota)\bigr)\subseteq F(\iota)$ and that
$y_\iota v=(v_\iota\otimes v{\restriction}_{E(\iota)})w_\iota$.
Furthermore, observe that for $a\in A_\iota$ and $\zeta\in E_\iota$,
$$ \bigl(v_\iota^*\sigma_\iota(a)v_\iota)\zeta
=v_\iota^*(a\otimes\zeta)=\theta_\iota(a)\zeta. $$
Hence for $a\in A_\iota$,
$$ \align
\theta(a)&=v^*\sigma(a)v
  =v^*\sigma_\iota(a)v
=v^*y_\iota^*\bigl(\sigma_\iota(a)\otimes1_{F(\iota)}\bigr)y_\iota v= \\
&=w_\iota^*\bigl(v_\iota^*\sigma_\iota(a)v_\iota
   \otimes(v{\restriction}_{E(\iota)})^*v{\restriction}_{E(\iota)}\bigr)w_\iota
  =w_\iota^*(\theta_\iota(a)\otimes1_{E(\iota)})w_\iota
  =\theta_\iota(a),
\endalign $$
and thus~(\thetathetai) holds.
Now to show that~(\thetaan) holds, consider
$a_j\in A_{\iota_j}\cap\ker\phi_{\iota_j}$ for some $\iota_j\in I$
($1\le j\le n$) with $\iota_j\ne\iota_{j+1}$.
It is easy to see that
$$ \aligned
\theta_{\iota_1}(a_1)\cdots\theta_{\iota_n}(a_n)\xi&=
\widehat{\theta_{\iota_1}(a_1)}\otdt\widehat{\theta_{\iota_n}(a_n)}= \\
&=v^*\bigl((a_1\otimes\xi_{\iota_1})\otdt(a_n\otimes\xi_{\iota_n})\bigr)
=\theta(a_1\cdots a_n)\xi. \endaligned
\tag{\thetaanxi} $$
Now consider an element $\zeta_1\otdt\zeta_p\in E$, where
$\zeta_j\in E\oup_{k_j}$ for some $k_j\in I$ with $k_j\ne k_{j+1}$.

Let $P_0:E\to\xi B$ be the projection and for $\ell\in\Nats$ let
$$ P_\ell:E\to\bigoplus\Sb
  \iota_1,\ldots,\iota_\ell\in I\\
  \iota_1\ne\iota_2,\ldots,\iota_{\ell-1}\ne\iota_\ell\endSb
  E_{\iota_1}\oup\otdt E_{\iota_\ell} $$
be the projection.
Taking adjoints and using~(\thetaanxi), we see that
$$ P_0\theta_{\iota_1}(a_1)\cdots\theta_{\iota_n}(a_n)(\zeta_1\otdt\zeta_p)
=P_0\theta(a_1\cdots a_n)(\zeta_1\otdt\zeta_p). $$
Now letting $\ell\in\Nats$ we will use standard techniques (see, for
example,~\cite{\DykemaHaagerupRordam} and~\cite{\DykemaZZExact}) to show that
$$ P_\ell\theta_{\iota_1}(a_1)\cdots\theta_{\iota_n}(a_n)(\zeta_1\otdt\zeta_p)
=P_\ell\theta(a_1\cdots a_n)(\zeta_1\otdt\zeta_p). \tag{\Plthetaan} $$
If $\ell>n+p$ or $\ell<|n-p|$ then it is clear that both sides of~(\Plthetaan)
are zero.
If $\ell=n+p$ then both sides of~(\Plthetaan) are zero unless $\iota_n\ne k_1$,
in which case a calculation similar to~(\thetaanxi) shows that~(\Plthetaan)
holds.
Let $Q\oup_\iota:E_\iota\to E_\iota\oup$ and
$R_\iota\oup:F_\iota\to F_\iota\oup$ be the projections, and note that
$R_\iota\oup v_\iota=v_\iota Q_\iota\oup$.
Consider when $n+p-\ell=1$.
Then both sides of~(\Plthetaan) are zero unless $\iota_n=k_1$, in which case
$$ \align
P_\ell\theta_{\iota_1}&(a_1)\cdots\theta_{\iota_n}(a_n)
  (\zeta_1\otdt\zeta_p)= \\
  \vspace{1ex}
&=\widehat{\theta_{\iota_1}(a_1)}\otdt\widehat{\theta_{\iota_{n-1}}(a_{n-1})}
  \otimes Q\oup_{\iota_n}(\theta_{\iota_n}(a_n)\zeta_1)
  \otimes\zeta_2\otdt\zeta_p \\
  \vspace{1ex} \allowdisplaybreak
&=v_\iota^*\bigl((a_1\otimes\xi_{\iota_1})\otdt
  (a_{n-1}\otimes\xi_{\iota_{n-1}})
  \otimes R_{\iota_n}\oup(a_n\otimes\zeta_1)
  \otimes(1\otimes\zeta_2)\otdt(1\otimes\zeta_p)\bigr) \\
  \vspace{1ex}
&=P_\ell\theta(a_1\cdots a_n)(\zeta_1\otdt\zeta_p) \endalign $$
If $p+n-\ell=2r+1$ for $r\in\{1,2,\ldots,\min(p,n)-2\}$ then both sides
of~(\Plthetaan) are zero unless
$\iota_n=k_1,\,\iota_{n-1}=k_2,\ldots,\iota_{n-r+1}=k_r$, in which case
$$ \align
P_\ell&\theta_{\iota_1}(a_1)\cdots\theta_{\iota_n}(a_n)
  (\zeta_1\otdt\zeta_p)= \\
  \vspace{2ex}
&=\topaligned
&\widehat{\theta_{\iota_1}(a_1)}\otdt\widehat{\theta_{\iota_{n-r-1}}(a_{n-r-1})}
    \otimes \\
   &\otimes Q_{\iota_{n-r}}\oup\bigl(\theta_{\iota_{n-r}}(a_{n-r})
    \langle\xi,\,\theta_{\iota_{n-r+1}}(a_{n-r+1})\cdots
    \theta_{\iota_n}(a_n)\zeta_1\otdt\zeta_r\rangle\zeta_{r+1}\bigr) \otimes \\
   &\otimes\zeta_{r+2}\otdt\zeta_p \endtopaligned \\
\vspace{2ex} \allowdisplaybreak
&=\topaligned
&\widehat{\theta_{\iota_1}(a_1)}\otdt\widehat{\theta_{\iota_{n-r-1}}(a_{n-r-1})}
    \otimes \\
   &\otimes Q_{\iota_{n-r}}\oup\bigl(\theta_{\iota_{n-r}}(a_{n-r})
    \langle\widehat{\theta_{\iota_n}(a^*_n)}\otdt
    \widehat{\theta_{\iota_{n-r+1}}(a_{n-r+1}^*)},\,
    \zeta_1\otdt\zeta_r\rangle\zeta_{r+1}
    \bigr)\otimes \\
   &\otimes\zeta_{r+2}\otdt\zeta_p \endtopaligned \\
\vspace{2ex} \allowdisplaybreak
&=\topaligned
   v^*\Bigl(&(a_1\otimes\xi_{\iota_1})\otdt(a_{n-r-1}\otimes\xi_{\iota_{n-r-1}})
    \otimes \\
   &\otimes \topaligned
        R_{\iota_{n-r}}\oup\Bigl(&\theta_{\iota_{n-r}}(a_{n-r})\;\cdot \\
        &\;\cdot
\bigl\langle\sigma_{\iota_n}(a_n^*)\cdots\sigma_{\iota_{n-r+1}}(a_{n-r+1}^*)\eta
        ,\,(1\otimes\zeta_1)\otdt(1\otimes\zeta_r)\bigr\rangle
        (1\otimes\zeta_{r+1})\Bigr)\otimes
        \endtopaligned \\
   &\otimes(1\otimes\zeta_{r+2})\otdt(1\otimes\zeta_p)\Bigr) \endtopaligned \\
\vspace{2ex} \allowdisplaybreak
&=\topaligned
   v^*\Bigl(&(a_1\otimes\xi_{\iota_1})\otdt(a_{n-r-1}\otimes\xi_{\iota_{n-r-1}})
    \otimes \\
   &\otimes \topaligned
        R_{\iota_{n-r}}\oup\Bigl(&\theta_{\iota_{n-r}}(a_{n-r})\;\cdot \\
          &\;\cdot
\bigl\langle\eta,\,\sigma_{\iota_{n-r+1}}(a_{n-r+1})\cdots\sigma_{\iota_n}(a_n)
        (1\otimes\zeta_1)\otdt(1\otimes\zeta_r)\bigr\rangle
        (1\otimes\zeta_{r+1})\Bigr)\otimes
        \endtopaligned \\
   &\otimes(1\otimes\zeta_{r+2})\otdt(1\otimes\zeta_p)\Bigr) \endtopaligned \\
\vspace{2ex}
&=P_\ell\theta(a_1\cdots a_n)(\zeta_1\otdt\zeta_p). \endalign $$
A similar calculation shows that~(\Plthetaan) holds also when
$n+p-\ell=2\min(p,n)-1$.

If $n+p-\ell=2r$ is even for $r\in\{1,2,\ldots,\min(p,n)-1\}$ then both sides
of~(\Plthetaan) are zero unless
$\iota_n=k_1,\,\iota_{n-1}=k_2,\ldots,\iota_{n-r+1}=k_r$ and
$\iota_{n-r}\ne k_{r+1}$, in which case
$$ \align
P_\ell\theta_{\iota_1}&(a_1)\cdots\theta_{\iota_n}(a_n)
   (\zeta_1\otdt\zeta_p)= \\
   \vspace{2ex}
&=\topaligned
   &\widehat{\theta_{\iota_1}(a_1)}\otdt\widehat{\theta_{\iota_{n-r}}(a_{n-r})}
    \otimes \\
   &\otimes\langle\xi,\,\theta_{\iota_{n-r+1}}(a_{n-r+1})\cdots
    \theta_{\iota_n}(a_n)\zeta_1\otdt\zeta_r\rangle\zeta_{r+1} \otimes \\
   &\otimes\zeta_{r+2}\otdt\zeta_p \endtopaligned \\
\vspace{2ex} \allowdisplaybreak
&=\topaligned
   &\widehat{\theta_{\iota_1}(a_1)}\otdt\widehat{\theta_{\iota_{n-r}}(a_{n-r})}
    \otimes \\
   &\otimes\langle\widehat{\theta_{\iota_n}(a^*_n)}\otdt
    \widehat{\theta_{\iota_{n-r+1}}(a_{n-r+1}^*)},\,
    \zeta_1\otdt\zeta_r\rangle\zeta_{r+1}
    \otimes \\
   &\otimes\zeta_{r+2}\otdt\zeta_p \endtopaligned \\
\vspace{2ex} \allowdisplaybreak
&=\topaligned
   v^*\Bigl(&(a_1\otimes\xi_{\iota_1})\otdt(a_{n-r}\otimes\xi_{\iota_{n-r}})
    \otimes \\
   &\otimes\bigl\langle \eta,\,\sigma_{\iota_{n-r+1}}(a_{n-r+1})\cdots
    \sigma_{\iota_n}(a_n)
    \bigl((1\otimes\zeta_1)\otdt(1\otimes\zeta_r)\bigr)\bigr\rangle
    (1\otimes\zeta_{r+1})
    \otimes \\
   &\otimes(1\otimes\zeta_{r+2})\otdt(1\otimes\zeta_p)\Bigr) \endtopaligned \\
   \vspace{2ex}
&=P_\ell\theta(a_1\cdots a_n)(\zeta_1\otdt\zeta_p). \endalign $$
Similar calculations show that~(\Plthetaan) holds also when
$p+n-\ell=2\min(p,n)$.
This finishes the proof of~(\thetaan), and of the theorem.
\QED

\heading \S\vNalgs.  Amalgamated free products of von Neumann
algebras. \endheading

This section contains results for amalgamated free products
of von Neumann algebras that are analogous to those for C$^*$--algebras found
in~\S\Embeddings{} and~\S\cpmaps.
The free product of von Neumann algebras with respect to given normal
states was defined by Voiculescu in~\cite{\VoiculescuZZSymmetries} and has been
much studied.
See also Ching's paper~\cite{\Ching}, where the free product of von Neumann
algebras with respect to normal faithful tracial states from a certain class
was first defined.
We begin this section by describing the ``folklore'' construction of 
amalgamated
free products of von Neumann algebras (with respect to normal conditional
expectations onto a von Neumann subalgebra).
We are grateful to the referee for pointing out this
simplification of the construction we originally gave.

\demo{Construction~\vNamFP}
Let $B$ be a von Neumann algebra contained as a unital von Neumann 
subalgebra of
von Neumann algebras $A_\iota$ ($\iota\in I$).
Suppose there are normal conditional expectations $\phi_\iota:A_\iota\to B$.
The amalgamated free product of von Neumann algebras, which we will denote
$$ (A,\phi)=\freeprodvni(A_\iota,\phi_\iota), $$
is constructed as follows.
Let
$$ (\Aeu,\varphi)=\freeprodi(A_\iota,\phi_\iota) $$
be the C$^*$--algebra reduced amalgamated free product.
Let $\psi$ be a normal state on $B$ with faithful GNS representation 
and consider the state $\psi\circ\varphi$ on $\Aeu$.
Let $\pi_{\psi\circ\varphi}$ be the GNS representation of $\Aeu$
associated to $\psi\circ\varphi$, let
$A=\pi_{\psi\circ\varphi}(\Aeu)''\subseteq\Leu(L^2(\Aeu,\psi\circ\varphi))$ 
and thereby regard $\Aeu$ as a
weakly dense subalgebra of $A$.
The Hilbert space projection $L^2(\Aeu,\psi\circ\varphi)\to 
L^2(B,\psi)$ gives rise to a normal
conditional expectation $\phi:A\to B$ whose restriction to $\Aeu$ is $\varphi$.
It remains to see that the pair $(A,\phi)$ is independent of the 
choice of $\psi$.
If $\psi'$ is any normal state on $B$ then the state 
$\psi'\circ\varphi$ of $\Aeu$
extends to the normal state $\psi'\circ\phi$ of $A$.
Hence $\pi_{\psi'\circ\varphi}(\Aeu)''=\pi_{\psi'\circ\phi}(A)$ is a 
quotient of $A$.
Thus $(A,\phi)$ is independent of the choice of $\psi$.
\enddemo

\demo{Remark \vNamFPrep}
In the above construction, we have $\Aeu\subseteq\Leu(E)$
for the Hilbert $B,B$--bimodule $E=L^2(\Aeu,\varphi)$.
The GNS Hilbert space $L^2(\Aeu,\psi\circ\varphi)$ is canonically isomorphic to
the internal tensor product
$E\otimes_{\pi_\psi}L^2(B,\psi)$, where 
$\pi_\psi:B\to\Leu(L^2(B,\psi))$ is the GNS representation of $\psi$,
and the representation $\pi_{\psi\circ\varphi}$ is
given by 
$\pi_{\psi\circ\varphi}(a)=a\otimes1\in\Leu(E\otimes_{\pi_\psi}L^2(B,\psi))$.
We will later use this picture in proofs.
\enddemo

Following Rieffel~\scite{\RieffelZZMoritaEqCW}{5.1},
if $A$ and $B$ are von Neumann algebras, if $E$ is a Hilbert $B$--module and if
$\theta:A\to\Leu(E)$ is a completely positive map, we say that $\theta$ is
{\it normal} if for every $\zeta_1,\zeta_2\in E$, the map
$A\ni a\mapsto \langle\zeta_1,\theta(a)\zeta_2\rangle\in B$ is normal.
This coincides with the usual notion of normality when $B=\Cpx$ (in which case
$E$ is a Hilbert space).
It is clear that if $B$ is a von Neumann subalgebra of a von Neumann algebra
$A$ having a normal conditional expectation $\phi:A\to B$ then the GNS
representation of $A$ as bounded adjointable operators on the Hilbert
$B$--module $L^2(A,\phi)$ is normal.

Part~(i) of the following straightforward lemma was proved in the case of a
$*$--homo\-morphism by Rieffel as part of~\scite{\RieffelZZMoritaEqCW}{5.2}.
\proclaim{Lemma \normality}
Let $A$ and $B$ be von Neumann algebras, let $E$ be a Hilbert $B$--module and
suppose that $\theta:A\to\Leu(E)$ is completely positive map.
Let $\Hil$ be a Hilbert space and let $\tau:B\to\Leu(\Hil)$ be a normal
$*$--representation.
Let $\theta\otimes1$ denote the completely positive map
$A\ni a\mapsto\pi(a)\otimes1\in\Leu(E\otimes_\tau\Hil)$
of $A$ into bounded operators on the Hilbert space $E\otimes_\tau\Hil$.
We have:
\roster
\item"(i)" if $\theta$ is normal then $\theta\otimes1$ is normal;
\item"(ii)" if $\theta\otimes1$ is normal and if $\tau$ is faithful then
$\theta$ is normal.
\endroster
\endproclaim
\demo{Proof}
If $\zeta_1,\zeta_2\in E$ and $v_1,v_2\in\Hil$ then
$$ \bigl\langle\zeta_1\otimes v_1,
  (\theta\otimes1)(x)(\zeta_2\otimes v_2)\bigr\rangle
  =\bigl\langle v_1,\tau\bigl(\langle\zeta_1,\theta(x)\zeta_2\rangle\bigr)
  v_2\bigr\rangle.
\tag{\evxev} $$
If $\theta$ is normal then~(\evxev) shows that $\theta\otimes1$ is
continuous from
the $\sigma(A,A_*)$ topology on $A$ to the weak--operator topology on
$\Leu(E\otimes_\tau\Hil)$, which implies $\theta\otimes1$ is normal and
proves~(i).

If $\theta\otimes1$ is normal then~(\evxev) shows that
$x\mapsto\tau(\langle\zeta_1,\theta(x)\zeta_2\rangle)$ is normal.
Assuming also $\tau$ is faithful, it follows that $\theta$ is normal.
\QED

The following application of Lemma~\normality(i) is
in~\scite{\RieffelZZMoritaEqCW}{5.2}.
\proclaim{Lemma \NormalInduced}
Let $B$ be a unital von Neumann subalgebra of a von Neumann algebra $A$ with a
normal conditional expectation $\Phi:A\to B$.
Let $\tau$ be a normal $*$--representation of $B$ on a Hilbert space $\Hil$.
Then the induced representation, $\tau{\induce}^A$, of $\tau$ to $A$ with
respect
to the conditional expectation $\Phi$ is normal.
\endproclaim
\demo{Proof}
By definition, and in the notation of Lemma~\normality,
$\tau{\induce}^A=\pi\otimes1:A\to\Leu(L^2(A,\Phi)\otimes_\tau\Hil)$, where
$\pi$ is the GNS representation of $A$ on $L^2(A,\Phi)$.
\QED

\proclaim{Lemma~\NormalTensor}
Let $A$, $B_1$ and $B_2$ be von Neumann algebras and let $E_j$ be a Hilbert
$B_j$--module ($j=1,2$).
If $\theta:A\to\Leu(E_1)$ and $\sigma:B_1\to\Leu(E_2)$ are normal completely
positive maps, then the completely positive map
$\theta\otimes1:A\to\Leu(E_1\otimes_\sigma E_2)$ is normal.
\endproclaim
\demo{Proof}
Let $\tau:B_2\to\Leu(\Hil)$ be a faithful normal $*$--representation of $B_2$
on a Hilbert space $\Hil$.
Applying Lemma~\normality{} in succession we find that
$\sigma\otimes1:B_1\to\Leu(E_2\otimes_\tau\Hil)$ is normal,
$\theta\otimes1\otimes1:A\to\Leu(E_1\otimes_\sigma E_2\otimes_\tau\Hil)$ is
normal, and thus $\theta\otimes1:A\to\Leu(E_1\otimes_\sigma E_2)$ is normal.
\QED

\proclaim{Lemma \vNCondExpA}
Let $B$ be a unital von Neumann algebra, let $I$ be a set and for every
$\iota\in I$
let $A_\iota$ be a unital von Neumann algebra containing a copy of $B$ as a
unital
von Neumann subalgebra and having a normal conditional expectation
$\phi_\iota:A_\iota\to B$ whose GNS representation is faithful.
Let
$$ (A,\phi)=\freeprodvni(A_\iota,\phi_\iota) $$
be the reduced amalgamated free product of von Neumann algebras.
Then for every $\iota_0\in I$, there is a normal conditional expectation
$\Phi_{\iota_0}:A\to A_{\iota_0}$ such that
$\Phi_{\iota_0}{\restriction}_{A_\iota}=\phi_\iota$ for every
$\iota\in I\backslash\{\iota_0\}$ and $\Phi_{\iota_0}(a_1a_2\cdots a_n)=0$
whenever $n\ge2$ and $a_j\in A_{\iota_j}\cap\ker\phi_{\iota_j}$ with
$\iota_1\ne\iota_2,\ldots,\iota_{n-1}\ne\iota_n$.
\endproclaim
\demo{Proof}
Let $\psi$ be a normal faithful state on $B$
and let $\tau=\pi_\psi:B\to\Leu(\Hil)$ be the associated GNS representation.
The construction of $A$ can be realized on the Hilbert space
$E\otimes_\tau\Hil$.
The projection $Q_{\iota_0}:E\to E_{\iota_0}$ from the proof of
Lemma~\CondExpA{} gives rise to the projection
$Q_{\iota_0}\otimes1_\Hil:E\otimes_\tau\Hil\to E_{\iota_0}\otimes_\tau\Hil$,
compression with which provides a normal positive linear map
$\Theta_{\iota_0}:A\to\Leu(E_{\iota_0}\otimes_\tau\Hil)$.
Let $\lambda_{\iota_0}:A_{\iota_0}\to\Leu(E_{\iota_0}\otimes_\tau\Hil)$ be the
GNS representation $A_{\iota_0}\hookrightarrow\Leu(E_{\iota_0})$ 
followed by the
inclusion
$\Leu(E_{\iota_0})\ni x\mapsto
x\otimes1_{\Hil}\in\Leu(E_{\iota_0}\otimes_\tau\Hil)$.
Then $\Theta_{\iota_0}$ maps a weakly dense $*$--subalgebra of $A$ into the
image of $\lambda_{\iota_0}$, hence maps all of $A$ there.
Let $\Phi_{\iota_0}=\lambda_{\iota_0}^{-1}\circ\Theta_{\iota_0}$.
The desired properties of $\Phi_{\iota_0}$ are easily verified.
\QED

Here is an embedding result, analogous to Theorem~\EmbThm,
for amalgamated free products of von Neumann algebras.
\proclaim{Theorem \vNEmbThm}
Let $B\subseteq\Bt$ be a (not necessarily unital) inclusion of von 
Neumann algebras.
Let $I$ be a set
and for each $\iota\in I$ suppose
$$ \matrix
1_{\At_\iota}\in& \Bt     & \subseteq & \At_\iota \\
                 & \cup &           & \cup  \\
1_{A_\iota}\in  & B       & \subseteq & A_\iota
\endmatrix $$
are inclusions of von Neumann algebras.
Suppose that $\phit_\iota:\At_\iota\to\Bt$ is a normal conditional 
expectation such that $\phit_\iota(A_\iota)\subseteq B$
and assume that $\phit_\iota$ and the restriction 
$\phit_\iota{\restriction}_{A_\iota}$ have faithful GNS 
representations,
for all $\iota\in I$.
Let
$$ \aligned
(\At,\phit)&=\freeprodvni(\At_\iota,\phit_\iota) \\
(A,\phi)&=\freeprodvni(A_\iota,\phit_\iota{\restriction}_{A_\iota})
\endaligned $$
be the amalgamated free products of von Neumann algebras.
Then there is a unique normal $*$--homo\-morph\-ism $\kappa:A\to\At$ such that
for every $\iota\in I$ the diagram
$$ \matrix
\At_\iota&\hookrightarrow&\At \\
\cup&&\phantom{\kappa}\uparrow\kappa \\
A_\iota&\hookrightarrow&A
\endmatrix \tag{\vNkappas} $$
commutes, where the horizontal arrows are the inclusions arising from 
the free product construction.
Moreover, $\kappa$ is necessarily injective.
\endproclaim
\demo{Proof}
This is very much like the proof of Theorem~\EmbThm, to which we refer in
detail.
Assume without loss of generality that $B$ is a unital subalgebra of $\Bt$.
We now insist that $\tau$ be a normal faithful representation of $\Bt$;
we must show that the algebra homomorphism
$\lambdat\circ\sigmat:\Afr\to\Leu(\Et\otimes_\tau\Weu)$ extends to a normal
representation of the von Neumann algebra $A$.
But $\Et\otimes_\tau\Weu$ is the direct sum of
$E\otimes_{\tau{\restriction}_B}\Weu$ and the various
$\Weut(\iota_1,\ldots,\iota_n)$.
The homomorphism $\lambdat\circ\sigmat$ restricted to
$E\otimes_{\tau{\restriction}_B}\Weu$
extends to the defining representation of $A$.
Let $n\ge1$ and let $\iota_1,\ldots,\iota_n\in I$ be such that
$\iota_j\ne\iota_{j+1}$;
we have the normal $*$-representation, $\mu$, of $A_{\iota_1}$ on the Hilbert
space $K_{\iota_1}\otimes_\Bt E_{\iota_2}\otdts\Bt 
E_{\iota_n}\otimes_\tau\Weu$, obtained from
the normal representation of $A_{\iota_1}$ in $\Leu(\Et_{\iota_1})$;
let $\mu{\induce}^A$ be the representation of the von Neumann algebra $A$ on
$\Weut(\iota_1,\ldots,\iota_n)$ induced from $\mu$ with respect to the normal
conditional expectation $\Phi_{\iota_1}:A\to A_{\iota_1}$ found in
Lemma~\vNCondExpA;
then $\lambdat\circ\sigmat$ restricted to $\Weut(\iota_1,\ldots,\iota_n)$
extends to the $*$--homomorphism 
$\mu{\induce}^A:A\to\Leu(\Weut(\iota_1,\ldots,\iota_n))$, which by
Lemma~\NormalInduced{} is normal.
\QED

Here is the construction, analogous to Theorem~\fpcpuBB, of free products of
certain completely positive maps in the von Neumann algebra setting.
\proclaim{Theorem \vNfpcpuBB}
Let $B$ be a von Neumann algebra, let $I$ be a set and for every $\iota\in I$
let $A_\iota$ and $D_\iota$ be von Neumann algebras containing copies of $B$
as unital von Neumann subalgebras and having normal conditional expectations
$\phi_\iota:A_\iota\to B$, respectively $\psi_\iota:D_\iota\to B$, whose GNS
representations are faithful.
Suppose that for each $\iota\in I$ there is a normal unital completely positive
map $\theta_\iota:A_\iota\to D_\iota$ that is also a $B$--$B$ bimodule map and
satisfies $\psi_\iota\circ\theta_\iota=\phi_\iota$.
Let
$$ \aligned
(A,\phi)&=\freeprodvni(A_\iota,\phi_\iota) \\
(D,\phi)&=\freeprodvni(D_\iota,\psi_\iota)
\endaligned \tag{\Advnfp} $$
be the reduced amalgamated free products of von Neumann algebras.
Then there is a normal unital completely positive map $\theta:A\to D$
satisfying
$$ \forall\iota\in I\qquad\theta{\restriction}_{A_\iota}=\theta_\iota
\tag{\vNthetathetai} $$
and
$$ \theta(a_1a_2\cdots a_n)=\theta(a_1)\theta(a_2)\cdots\theta(a_n)
\tag{\vNthetaan} $$
whenever $a_j\in A_{\iota_j}\cap\ker\phi_{\iota_j}$ and
$\iota_1\ne\iota_2$, $\iota_2\ne\iota_3,\,\ldots,\,\iota_{n-1}\ne\iota_n$.
\endproclaim
\demo{Proof}
Let
$$ \aligned
(\Aeu,\varphi)&=\freeprodi(A_\iota,\phi_\iota) \\
(\Deu,\rho)&=\freeprodi(D_\iota,\psi_\iota)
\endaligned $$
be the C$^*$--algebra reduced amalgamated free products.
Thus $A$ and $D$ are the closures in strong--operator topology of $\Aeu$
and respectively $\Deu$ in the appropriate representations.
We need only show that the unital completely positive map 
$\theta:\Aeu\to\Deu$ found
in Theorem~\fpcpuBB{} extends to a normal completely positive map
$\thetabar:A\to D$.

Consider, from the proof of Theorem~\fpcpuBB, the Hilbert $B$--modules
$(E,\xi)=\freeprodi(E_\iota,\xi_\iota)$, and
$(F,\eta)=\freeprodi(F_\iota,\eta_\iota)$, the $*$--homomorphism
$\sigma:\Aeu\to\Leu(F)$ and the bounded operator $v\in\Leu(E,F)$;
denote by $i_\Aeu$ the GNS representation of $\Aeu$ on $L^2(\Aeu,\varphi)$.
Recall that $\sigma$ is a free product of embeddings
$\sigma_\iota:A_\iota\to\Leu(F_\iota)$.
{}From the proof of Theorem~\EmbThm, and
letting $\tau=\pi_\mu:B\to\Leu(\Veu)$ be the GNS representation of a 
normal faithful state $\mu$ of $B$,
we see that the representation $\sigma\otimes1$ of $\Aeu$ on
$F\otimes_\tau\Veu$ given by $a\mapsto \sigma(a)\otimes1$ splits as a direct
sum,
$\sigma\otimes1=\bigoplus_{\lambda\in\Lambda}
(\sigma\otimes1){\restriction}_{\Weu_\lambda}$,
where each summand $(\sigma\otimes1){\restriction}_{\Weu_\lambda}$ is either a
copy of 
$i_\Aeu\otimes1:\Aeu\to\Leu(L^2(\Aeu,\varphi)\otimes_\tau\Veu)$ or is 
the induced
representation $\nu{\induce}^\Aeu$ of a representation $\nu$ of some 
$A_\iota$ on
a Hilbert space, where $\nu$ is the restriction to an invariant subspace of
the representation
$\sigma_\iota\otimes1:A_\iota\to\Leu(F_\iota\otimes_\tau\Veu)$.
The representation $i_\Aeu\otimes1$ extends to a normal $*$--representation of
$A$ as seen in Construction~\vNamFP.
Using Lemma~\normality{} we see that $\sigma_\iota\otimes1$ is normal;
hence $\nu$ is normal and by Lemma~\NormalInduced,
$\nu{\induce}^\Aeu$ extends to
a normal $*$--representation of $A$.
Hence $\sigma\otimes1$ extends to a normal $*$--representation of $A$,
which we will denote by $\sigmabar:A\to\Leu(F\otimes_\tau\Veu)$.

The isometry $v\in\Leu(E,F)$ gives rise to an isometry
$v\otimes1:E\otimes_\tau\Veu\to F\otimes_\tau\Veu$.
Letting $i_D:\Deu\to\Leu(E)$ be the defining representation, by
Remark~\vNamFPrep{}
the weak closure of the
image of $i_D\otimes1:\Deu\to\Leu(E\otimes_\tau\Veu)$ is
the von Neumann algebra $D$.
Consider the normal unital completely positive map
$\thetabar:A\to\Leu(E\otimes_\tau\Veu)$ given by
$\thetabar(x)=(v\otimes1)^*\sigmabar(x)(v\otimes1)$.
If $a\in\Aeu$ then $\thetabar(a)=i_D(\theta(a))\otimes1$.
So $\thetabar$ extends the map $\theta:\Aeu\to\Deu$;
hence $\thetabar(A)\subseteq D$.
\QED

\enddocument